\magnification=1200


\def\item{\vskip1.3pt\hang\textindent}


\tolerance=300 \pretolerance=200 \hfuzz=1pt \vfuzz=1pt

\hoffset 0cm            
\hsize=5.8 true in \vsize=9.5 true in

\def\rightheadline{\hfil\smc\lastname\hfil\tenbf\folio}
\def\leftheadline{\tenbf\folio\hfil\smc\lastname\hfil}
\headline={\ifodd\pageno\rightheadline\else\leftheadline\fi}
\newdimen\dimenone
\def\checkleftspace#1#2#3#4#5{
 \dimenone=\pagetotal
 \advance\dimenone by -\pageshrink   
 \ifdim\dimenone>\pagegoal          
   \else\dimenone=\pagetotal
        \advance\dimenone by \pagestretch
        \ifdim\dimenone<\pagegoal
          \dimenone=\pagetotal
          \advance\dimenone by#1         
          \setbox0=\vbox{#2\parskip=0pt                
                       \hyphenpenalty=10000
                       \rightskip=0pt plus 5em
                       \noindent#3 \vskip#4}    
        \advance\dimenone by\ht0
        \advance\dimenone by 3\baselineskip
        \ifdim\dimenone>\pagegoal\vfill\eject\fi
          \else\eject\fi\fi}

\parindent=35pt
\mathsurround=1pt
\parskip=1pt plus .25pt minus .25pt
\normallineskiplimit=.99pt

\mathchardef\emptyset="001F 

\def\Int{\mathop{\rm int}\nolimits}
%



\def\1{{\bf1}}\def\0{{\bf0}}

\def\({\bigl(}  \def\){\bigr)}
\def\<{\mathopen{\langle}}\def\>{\mathclose{\rangle}}

\def\Z{{\mathchoice{{\hbox{$\rm Z\hskip 0.26em\llap{\rm Z}$}}}%
{{\hbox{$\rm Z\hskip 0.26em\llap{\rm Z}$}}}%
{{\hbox{$\scriptstyle\rm Z\hskip 0.31em\llap{$\scriptstyle\rm Z$}$}}}{{%
\hbox{$\scriptscriptstyle\rm
Z$\hskip0.18em\llap{$\scriptscriptstyle\rm Z$}}}}}}

\def\F{{\mathchoice{\hbox{$\rm I\hskip-0.14em F$}}%
{\hbox{$\rm I\hskip-0.14em F$}}%
{\hbox{$\scriptstyle\rm I\hskip-0.14em F$}}%
{\hbox{$\scriptscriptstyle\rm I\hskip-0.10em F$}}}}

\def\R{{\mathchoice{\hbox{$\rm I\hskip-0.14em R$}}%
{\hbox{$\rm I\hskip-0.14em R$}}%
{\hbox{$\scriptstyle\rm I\hskip-0.14em R$}}%
{\hbox{$\scriptscriptstyle\rm I\hskip-0.10em R$}}}}

\def\C{{\mathchoice%
{\hbox{$\rm C\hskip-0.47em\hbox{%
\vrule height 0.58em width 0.06em depth-0.035em}$}\;}%
{\hbox{$\rm C\hskip-0.47em\hbox{%
\vrule height 0.58em width 0.06em depth-0.035em}$}\;}%
{\hbox{$\scriptstyle\rm C\hskip-0.46em\hbox{%
$\scriptstyle\vrule height 0.365em width 0.05em depth-0.025em$}$}\;}
{\hbox{$\scriptscriptstyle\rm C\hskip-0.41em\hbox{
$\scriptscriptstyle\vrule height 0.285em width 0.04em
depth-0.018em$}$}\;}}}

\def\.{{\cdot}}
\def\|{\Vert}
\def\ssk{\smallskip}
\def\msk{\medskip}
\def\bsk{\bigskip}
\def\giantskip{\vskip2\bigskipamount}

\def\giantbreak{\par \ifdim\lastskip<2\bigskipamount \removelastskip
         \penalty-400 \giantskip\fi}

\def\nin{\noindent}
\def\cen{\centerline}
\def\pagebreak{\vskip 0pt plus 0.0001fil\break}
\def\linebreak{\break}

\def\epsilon{\varepsilon}

\font\ninerm=cmr9 \font\eightrm=cmr8 \font\sixrm=cmr6
 \font\eightbf=cmbx8 \font\sixbf=cmbx6
 \font\eighti=cmmi8 \font\sixi=cmmi6
\font\ninesy=cmsy9 \font\eightsy=cmsy8 \font\sixsy=cmsy6
 \font\eightit=cmti8 
 \font\eightsl=cmsl8 
\font\eighttt=cmtt8
\font\bfone=cmbx10 scaled\magstep1 
\font\smc=cmcsc10 
 
scaled\magstep1 \font\small=cmcsc8

\def\no #1. {\bigbreak\vskip-\parskip\noindent\bf #1. \quad\rm}

\def\Proposition #1. {\checkleftspace{0pt}{\bf}{Theorem}{0pt}{}
\bigbreak\vskip-\parskip\noindent{\bf Proposition #1.} \quad\it}

\def\Theorem #1. {\checkleftspace{0pt}{\bf}{Theorem}{0pt}{}
\bigbreak\vskip-\parskip\noindent{\bf  Theorem #1.} \quad\it}
\def\Corollary #1. {\checkleftspace{0pt}{\bf}{Theorem}{0pt}{}
\bigbreak\vskip-\parskip\nin{\bf Corollary #1.} \quad\it}
\def\Lemma #1. {\checkleftspace{0pt}{\bf}{Theorem}{0pt}{}
\bigbreak\vskip-\parskip\noindent{\bf  Lemma #1.}\quad\it}

\def\Definition #1. {\checkleftspace{0pt}{\bf}{Theorem}{0pt}{}
\rm\bigbreak\vskip-\parskip\noindent{\bf Definition #1.} \quad}

\def\Remark #1. {\checkleftspace{0pt}{\bf}{Theorem}{0pt}{}
\rm\bigbreak\vskip-\parskip\noindent{\bf Remark #1.}\quad}

\def\Exercise #1. {\checkleftspace{0pt}{\bf}{Theorem}{0pt}{}
\rm\bigbreak\vskip-\parskip\noindent{\bf Exercise #1.} \quad}

\def\Example #1. {\checkleftspace{0pt}{\bf}{Theorem}{0pt}{}
\rm\bigbreak\vskip-\parskip\noindent{\bf Example #1.}\quad}
\def\Examples #1. {\checkleftspace{0pt}{\bf}{Theorem}{0pt}
\rm\bigbreak\vskip-\parskip\noindent{\bf Examples #1.}\quad}

\newcount\problemnumb \problemnumb=0
\def\Problem{\global\advance\problemnumb by 1\bigbreak\vskip-\parskip\noindent
{\bf Problem \the\problemnumb.}\quad\rm }

\def\Proof#1.{\rm\par\ifdim\lastskip<\bigskipamount\removelastskip\fi\smallskip
            \noindent {\bf Proof.}\quad}

\nopagenumbers

\def\author{}
\def\lastname{}
\def\thanks#1{\footnote*{\eightrm#1}}
\def\title{}

\def\nonumbers{\def\leftheadline{\hfil} \def\rightheadline{\hfil}}

\def\lastname{}
\def\h{{\textstyle{1\over2}}}

\def\he{{1\over2}}

\def\ep{\epsilon}

\def\text{\textstyle}
\def\disp{\displaystyle}
\def\d{{\,\rm d}}

\def\and{{\rm and }}

\expandafter\edef\csname amssym.def\endcsname{%
       \catcode`\noexpand\@=\the\catcode`\@\space}
\catcode`\@=11
\def\undefine#1{\let#1\undefined}
\def\newsymbol#1#2#3#4#5{\let\next@\relax
 \ifnum#2=\@ne\let\next@\msafam@\else
 \ifnum#2=\tw@\let\next@\msbfam@\fi\fi
 \mathchardef#1="#3\next@#4#5}
\def\mathhexbox@#1#2#3{\relax
 \ifmmode\mathpalette{}{\m@th\mathchar"#1#2#3}%
 \else\leavevmode\hbox{$\m@th\mathchar"#1#2#3$}\fi}
\def\hexnumber@#1{\ifcase#1 0\or 1\or 2\or 3\or 4\or 5\or 6\or 7\or 8\or
 9\or A\or B\or C\or D\or E\or F\fi}


\font\tenmsb=msbm10 \font\sevenmsb=msbm7 \font\fivemsb=msbm5
\newfam\msbfam
\textfont\msbfam=\tenmsb\scriptfont\msbfam=\sevenmsb
\scriptscriptfont\msbfam=\fivemsb \edef\msbfam@{\hexnumber@\msbfam}
\def\Bbb#1{{\fam\msbfam\relax#1}}

\font\teneufm=eufm10 \font\seveneufm=eufm7 \font\fiveeufm=eufm5
\newfam\eufmfam
\textfont\eufmfam=\teneufm \scriptfont\eufmfam=\seveneufm
\scriptscriptfont\eufmfam=\fiveeufm

\catcode`@=11 

\expandafter\edef\csname amssym.def\endcsname{%
       \catcode`\noexpand\@=\the\catcode`\@\space}
\font\eightmsb=msbm8 \font\sixmsb=msbm6 \font\fivemsb=msbm5
\font\eighteufm=eufm8 \font\sixeufm=eufm6 \font\fiveeufm=eufm5
\newskip\ttglue
\def\eightpoint{\def\rm{\fam0\eightrm}%
  \textfont0=\eightrm \scriptfont0=\sixrm \scriptscriptfont0=\fiverm
  \textfont1=\eighti \scriptfont1=\sixi \scriptscriptfont1=\fivei
  \textfont2=\eightsy \scriptfont2=\sixsy \scriptscriptfont2=\fivesy
  \textfont3=\tenex \scriptfont3=\tenex \scriptscriptfont3=\tenex
\textfont\eufmfam=\eighteufm \scriptfont\eufmfam=\sixeufm
\scriptscriptfont\eufmfam=\fiveeufm \textfont\msbfam=\eightmsb
\scriptfont\msbfam=\sixmsb \scriptscriptfont\msbfam=\fivemsb
  \def\it{\fam\itfam\eightit}%
  \textfont\itfam=\eightit
  \def\sl{\fam\slfam\eightsl}%
  \textfont\slfam=\eightsl
  \def\bf{\fam\bffam\eightbf}%
  \textfont\bffam=\eightbf \scriptfont\bffam=\sixbf
   \scriptscriptfont\bffam=\fivebf
  \def\tt{\fam\ttfam\eighttt}%
  \textfont\ttfam=\eighttt
  \tt \ttglue=.5em plus.25em minus.15em
  \normalbaselineskip=9pt
  \def\MF{{\manual opqr}\-{\manual stuq}}%
  \let\big=\eightbig
  \setbox\strutbox=\hbox{\vrule height7pt depth2pt width\z@}%
  \normalbaselines\rm}
\def\eightbig#1{{\hbox{$\textfont0=\ninerm\textfont2=\ninesy
  \left#1\vbox to6.5pt{}\right.\n@space$}}}


\csname amssym.def\endcsname


\def\al{\alpha}
\def\be{\beta}

\def\({\left(}
\def\){\right)}
\def\for{\qquad \hbox{for}\ }
\def\eq{\eqalign}

\def\O#1{O\(#1\)}
\def\abs#1{\left| #1 \right|}

\def\norm#1{\left\Vert #1 \right\Vert}

\def\klein{\eightpoint \def\smc{\small} \baselineskip=9pt}

\def\fn#1#2{{\parindent=0.7true cm
\footnote{$^{(#1)}$}{{\klein  #2}}}}

\font\boldmas=msbm10                  
\def\Bbb#1{\hbox{\boldmas #1}}        
\def\Z{{\Bbb Z}}                        

\def\R{{\Bbb R}}
\def\F{{\Bbb F}}

\def\C{{\Bbb C}}


\font\eightrm=cmr8 \long\def\fussnote#1#2{{\baselineskip=9pt
\setbox\strutbox=\hbox{\vrule height 7pt depth 2pt width 0pt}%
\eightrm \footnote{#1}{#2}}}
\font\boldmasi=msbm10 scaled 700      
\def\Bbbi#1{\hbox{\boldmasi #1}}      
\font\boldmas=msbm10                  
\def\Bbb#1{\hbox{\boldmas #1}}        
\def\Zi{{\Bbbi Z}}                      
\def\Pi{{\Bbbi P}}                      
\def\Ri{{\Bbbi R}}



\def\dint #1 {
\quad  \setbox0=\hbox{$\disp\int\!\!\!\int$}
  \setbox1=\hbox{$\!\!\!_{#1}$}
  \vtop{\hsize=\wd1\centerline{\copy0}\copy1} \quad}

\def\drint #1 {
\qquad  \setbox0=\hbox{$\disp\int\!\!\!\int\!\!\!\int$}
  \setbox1=\hbox{$\!\!\!_{#1}$}
  \vtop{\hsize=\wd1\centerline{\copy0}\copy1}\qquad}

\def\frac#1#2{{#1\over #2}}

\def\date{\the\day.~\the\month.~\the\year}

\def\mod{\,{\rm mod}\,}
\def\klein{\eightpoint \def\smc{\small} }

\def\frac#1#2{{#1\over#2}}
\def\Int{\int\limits}

\def\vol{{\rm vol}}

\nonumbers

\hsize=16.4true cm     \vsize=23.3true cm

\parindent=0cm

\def\eqno{\leqno}
\def\l{\ell}
\def\DS{\sum_{1\le h\le U}\ \sum_{1\le m\le u}}
\def\Sh{\sum_{h=1}^{[U]}}
\def\alh{\al_{h,[U]}}  \def\beh{\be_{h,[U]}} \def\gah{\gamma_{h,[U]}}
\def\Mj{{\cal M}_j}
\def\J{{\cal J}}
\def\D{{\cal D}}
\def\E{{\cal E}}
\def\F{{\cal F}}
\def\odd{\ {\rm odd}}
\def\b#1{{\bf #1}}
\def\rr{r_1\cdot\dots\cdot r_\l}

\vbox{\vskip 1.5true cm}

\footnote{}{\klein{\it Mathematics Subject Classification }
(2000): 11N37, 35P20, 58J50, 11P21.\par }

\cen{{\bfone A lower bound for the error term in Weyl's law}} \msk
\cen{{\bfone for certain Heisenberg manifolds}}\bsk \cen{{\bf
 Werner Georg Nowak }\fn{*}{The author gratefully
acknowledges support from the Austrian Science Fund (FWF) under
project Nr.~P20847-N18.} {\bf(Vienna)} }

\vbox{\vskip 1.2true cm}

{\klein{\bf Abstract. } This article provides an Omega-result for
the remainder term in Weyl's law for the spectral counting
function of certain rational $(2\l+1)$-dimensional Heisenberg
manifolds.}

\vbox{\vskip 1true cm}

{\bf Introduction. } For $M$ a closed $n$-dimensional Riemannian
manifold with a metric $g$ and Laplace-Beltrami operator $\Delta$,
let $N(t)$ denote the spectral counting function $$
N(t):=\sum_{\lambda{\ \rm eigenvalue\ of\ }\Delta\atop \lambda\le t}
d(\lambda) $$ where $d(\lambda)$ is the dimension of the eigenspace
corresponding to $\lambda$, and $t$ is a large real variable. Then a
deep and very general theorem due to L.~H\"{o}rmander [6] tells us that
$$ N(t) = {\vol(M)\over(4\pi)^{n/2}\,\Gamma(\h n+1)}\,t^{n/2} +
\O{t^{(n-1)/2}} \eqno(1.1) $$ ("Weyl's law"), and that the error
term in general cannot be improved. Nevertheless, it is of interest
to study the order of magnitude and the asymptotic behavior of the
remainder $ R(t) = N(t)- {\vol(M)\over(4\pi)^{n/2}\,\Gamma(\he
n+1)}\,t^{n/2}$ for special manifolds $M$. \ssk The most classic
example, namely the case that $M=\R^n/\Z^n$, the $n$-dimensional
torus, is (equivalent to) a central problem in the theory of lattice
points in large domains, namely to provide asymptotic results for
the number $A_n(x)$ of integer points in an origin-centered
$n$-dimensional ball of radius $x$, for any dimension $n\ge2$. There
exists a vast literature on this particular subject: We only refer
to the works of Huxley [7], [8], Hafner [4], and Soundararajan [17]
for the planar case, for the papers by Chamizo \& Iwaniec [1],
Heath-Brown [5], and Tsang [18] for dimension $n=3$, and to the
monographs of Walfisz [20], and Kr\"{a}tzel [14], [15], as well as to
the recent, quite comprehensive, survey article [9]. \ssk In fact,
for $M=\R^n/\Z^n$, we see that\fn{1}{Bold face letters will denote
throughout elements of $\Ri^n$, resp., $\Zi^n$, which may be viewed
also as $(1\times n)$-matrices ("row vectors") where applicable.
Further, $|\cdot|$ stands for the Euclidean norm.} $\{\b u\mapsto
e(\b m\cdot\b u ):\ \b m\in\Z^n\}$ is a basis for the eigenfunctions
of the Laplace operator $\Delta=-\sum_{j=1}^n\partial_{jj}$, acting
on functions from $\R^n/\Z^n$ into $\C$. The corresponding
eigenvalues are $4\pi^2|\b m|^2$, $\b m\in\Z^n$. For any integer
$k\ge0$, let as usual $r_n(k)$ denote the number of ways to write
$k$ as the sum of $n$ squares. Then, for each $k$ with $r_n(k)>0$,
$4\pi^2 k$ is an eigenvalue of $\Delta$ whose eigenspace consists of
all functions $$ \b u \mapsto \sum_{|\b m |^2=k}c(\b m)\,e(\b
m\cdot\b u)\,, $$ where $c(\b m)$ are any complex coefficients. Its
dimension obviously equals $r_n(k)$, hence
$$ N(t) = \sum_{k\ge0:\ 4\pi^2 k\le t} r_n(k) =
A_n\({\sqrt{t}\over2\pi}\)\,. $$ \bsk \msk

{\bf 2. Heisenberg manifolds. } In recent times, presumably
motivated by quite different areas like quantum physics and the
abstract theory of PDE's, a lot of work has been done on another
special case, namely that of so-called Heisenberg manifolds. To
recall basics and to fix notions, let $\l\ge1$ be a given integer,
and
$$ \gamma(\b x,\b y,z) = \pmatrix{1 &\b x&z\cr ^t\b o_\l&I_\l& ^t\b
y\cr 0&\b o_\l&1\cr}\,,
$$ where $\b x, \b y \in \R^\l$, $z\in\R$,
$\b o_\l =(0,\dots,0)\in\R^\l$, $I_\l$ is the $(\l\times\l)$-unit
matrix, and $^t\cdot$ denotes transposition. Then the
$(2\l+1)$-dimensional Heisenberg group $H_\l$ is defined by $$ H_\l
= \{\gamma(\b x,\b y,z):\ \b x, \b y \in \R^\l,\ z\in\R\ \}\,,
\eqno(2.1) $$ with the usual matrix product. Further, for any
$\l$-tuple $\b r=(r_1,\dots,r_\l)\in\Z_+^\l$ with the property that
$r_j\mid r_{j+1}$ for all $j=1,\dots,\l-1$, we put $\b r * \Z^\l:=
r_1\Z\times\dots\times r_\l\Z $ and define $$ \Gamma_{\b r} =
\{\gamma(\b x,\b y,z):\ \b x \in \b r*\Z^\l, \b y \in\Z^\l , z\in\Z\
\}\,. \eqno(2.2)$$ $\Gamma_{\b r}$ is a uniform discrete subgroup of
$H_\l$, i.e., the {\it Heisenberg manifold } $H_\l/\Gamma_{\b r}$ is
compact. Fortunately, according to a deep work by Gordon and Wilson
[2], Theorem 2.4, this seemingly quite special choice of $\Gamma_{\b
r}$ is in fact fairly general. The subgroups $\Gamma_{\b r}$
classify {\it all } uniform discrete subgroups of $H_\l$ up to
automorphisms: For every uniform discrete subgroup $\Gamma$ of
$H_\l$ there exists a unique $\l$-tuple $\b r$ and an automorphism
of $H_\l$ which maps $\Gamma$ to $\Gamma_{\b r}$. \ssk However, to
get a "rational" or "arithmetic" Heisenberg manifold - borrowing an
expression due to Petridis \& Toth [16] - we have to make a quite
particular choice of the metric involved\fn{2}{Compare the
discussion below concerning the bound (3.3) valid for "almost all"
metrics $g$.}. Following the example of [16], Theorem 1.1, and also
Zhai [21], we pick
$$ g_\l = \pmatrix{I_{2\l}& ^t\b o_{2\l}\cr \b o_{2\l}&2\pi}\,. \eqno(2.3)$$
The spectrum of the Laplace-Beltrami operator on $H_\l/\Gamma_{\b
r}$ has been analyzed in Gordon and Wilson [2], p.~259, and also
in Khosravi and Petridis [12], p.~3564. It consists of two
different classes ${\cal S}_I$ and ${\cal S}_{II}$, where ${\cal
S}_{I}$ is the spectrum of the Laplacian on the $2\l$-dimensional
torus, and
$$ {\cal S}_{II} = \{2\pi\(n_0^2+ n_0(2n_1+\l)\)\,:\
n_0\in\Z^+,\ n_1\in\Z_0^+ \ \}\,, $$ with multiplicities ( =
dimensions of corresponding eigenspaces)
$2n_0^\l\,\rr\,{n_1+\l-1\choose\l-1}$. \bsk\msk

{\bf 3.~Statement of problem and results. } In this article, we
shall be concerned with the size of the error term in (1.1) for the
special case that $M=(H_\l/\Gamma_{\b r}, g_\l)$ as described above,
i.e., with the asymptotic behavior of $$ R(t) = N(t)-
{\vol(M)\over(4\pi)^{\l+1/2}\,\Gamma(\l+{3\over2})}\,t^{\l+1/2}=
N(t) - {\rr\over2^{2\l+1/2}\pi^{\l}\,
\Gamma(\l+{3\over2})}\,t^{\l+1/2}\,. \eqno(3.1) $$ For $\l=1$,
Petridis and Toth [16] proved that $R(t)\ll t^{5/6}\log t$. They
were the first to realize that the question is related to a certain
planar lattice point problem which can be dealt with usual tools for
the estimation of fractional part sums. This result was sharpened
and generalized to arbitrary $\l\ge1$ by Khosravi and Petridis [12]
who obtained $R(t)\ll t^{\l-7/41}$. Zhai [21] noticed that Huxley's
deep method [7], [8] can be used to derive, for any $\l\ge1$,
$$ R(t) \ll t^{\l-77/416}(\log t)^{26957/8320}\,. \eqno(3.2)$$ In
fact, the difficulty in these estimations comes from the
"rational" nature of the metric $g_\l$. As Khosravi and Petridis
[12] showed, for "almost all" metrics $g$ the much sharper bound
$$ R_g(t)\ll_g\ t^{\l-1/4}\,\log t \eqno(3.3) $$ holds true.
Returning to the rational case (2.3), a result of Khosravi [11]
and Khosravi \& Toth [13] tells us that $$ \Int_0^T (R(t))^2 \d t
= C_\l\,T^{2\l+1/2} + \O{T^{2\l+1/4+\ep}} \eqno(3.4) $$ where
$C_\l>0$ is an explicit constant. A recent paper of Zhai [21]
provides estimates and asymptotics for higher power moments of
$R(t)$. In fact, (3.3) and (3.4) may suggest the conjecture that
$$ R(t) \ll t^{\l-1/4+\ep} \eqno(3.5) $$ for every $\ep>0$. The
situation has a good deal in common with the Dirichlet divisor and
the Gaussian circle problems. \ssk The objective of the present
article is to provide a lower bound which shows that the $\ep$ in
(3.5) cannot be removed. We shall prove that $$ R(t)=
\Omega\(t^{\l-1/4}\,(\log t)^{1/4}\)\,. $$ Together with (3.4), we
may say that "$R(t) \ll t^{\l-1/4}$ {\it in mean-square, } with an
unbounded sequence of exceptionally large values $t$". \ssk
Unfortunately, we have to impose the restriction that $\l$ is an
{\it even } integer. We will comment on this condition at the end of
the paper. \bsk

\bsk {\bf Theorem. } { \it For a fixed even positive integer $\l$,
let $(H_\l/\Gamma_{\b r},g_\ell)$ be a rational
$(2\l+1)$-dimensional Heisenberg manifold with metric $g_\l$, as
described above. Then the error term $R(t)$ for the associated
spectral counting function, defined in $(3.1)$, satisfies
$$ \limsup_{t\to\infty}{R(t)\over t^{\l-1/4}\,(\log t)^{1/4}} > 0\,. $$ } \bsk\msk

{\bf 4.~Some Lemmas.}\bsk

{\bf Lemma 1. } (Vaaler's approximation of fractional parts by
trigonometric polynomials.) { \it For arbitrary $w\in\R$ and
$H\in\Z^+$, let $\psi(w):= w-[w]-\h$, $$ \Sigma_H(w):=\sum_{h=1}^H
{\al_{h,H}}\,\sin(2\pi hw)\,,\qquad \Sigma_H^*(w):=\sum_{h=1}^H
{\be_{h,H}}\,\cos(2\pi hw)\ + {1\over2H+2}\,, $$ where, for
$h=1,\dots,H$,
$$ \al_{h,H}:={1\over\pi h}\,\rho\({h\over H+1}\)\,,\qquad \be_{h,H}:={1\over H+1}\(1-{h\over
H+1}\)\,, $$ and $$ \rho(\xi)= \pi\xi(1-\xi)\cot(\pi\xi)+\xi
\qquad\qquad (0<\xi<1)\,. $$ Then the following inequality holds
true:
$$ \abs{\psi(w)+\Sigma_H(w)} \le \Sigma_H^*(w)\,. $$ } \bsk

{\bf Proof. } This is one of the main results in Vaaler [19]. A
very well readable exposition can also be found in the monograph
by Graham and Kolesnik [3].

 \bsk\msk

{\bf Lemma 2. } { \it Let $F\in C^4[A, B]$,  $G\in C^2[A, B]$, and
suppose that, for positive parameters $X, Y, Z$, we have $1\ll
B-A\ll X$ and
$$ F^{(j)}\ll X^{2-j} Y^{-1} \for j=2,3,4, \
\abs{F''}\ge c_0 Y^{-1}\,,\quad G^{(j)}\ll X^{-j} Z \for j=0,1,2,
$$ throughout the interval $[A,B]$, with some constant $c_0>0$.
Let $\J'$ denote the image of $]A,B]$ under $F'$, and $F^*$ the
inverse function of $F'$. Then, with $e(w)=e^{2\pi iw}$ as usual,
$$ \eq{\sum_{A<m\le B} G(m)\,e(F(m)) =&\ e\({{\rm sgn}(F'')\over8}\) \sum_{k\in\J'}{G(F^*(k))
\over\sqrt{\abs{F''(F^*(k))}}}\,e\(F(F^*(k))-kF^*(k)\)  + \cr & +
\O{Z\(\sqrt{Y}+\log(2+\,{\rm length}(\J'))\)}\,. \cr }
$$ }\bsk

{\bf Proof. } Transformation formulas of this kind are quite
common, though often with worse error terms. This very sharp
version can be found as f.~(8.47) in the recent monograph [10] of
H.~Iwaniec and E.~Kowalski.

\bsk\msk

{\bf Lemma 3. } {\it For a real parameter $T\ge1$, let $\F_T$
denote the Fej\'er kernel $$ \F_T(v) = T\({\sin(\pi Tv)\over\pi Tv
}\)^2\,. $$ Then for arbitrary real numbers $Q>0$ and $\delta$, it
follows that $$ \Int_{-1}^1 \F_T(v) e(Qv+\delta)\d v =
\max\(1-{Q\over T},0\)e(\delta) + \O{1\over Q}\,, $$ where the
$O$-constant is independent of $T$ and $\delta$.}\bsk

{\bf Proof. } This useful result is due to Hafner [4]. It follows
from the classic Fourier transform formula $$ \Int_\Ri \F_T(v)
e(Qv) \d v = \Int_\Ri \({\sin(\pi v)\over\pi v }\)^2\ e\({Q\over
T}v\)\d v = \max\(1-{Q\over T},0\)\,. $$ Since $\F_T(\pm1)\ll
T^{-1}$ and $\F_T'(v)\ll v^{-2}$ for $|v|\ge1$, uniformly in
$T\ge1$, integration by parts readily shows that the intervals
$]-\infty,-1]$ and $[1,\infty[$ contribute only $\O{Q^{-1}}$.
\bsk\bsk

{\bf 5.~Proof of the Theorem. }\quad We start from Lemma 3.1 in
Zhai [21] which approximates the error term involved by a
fractional part sum. Let $U$ be a large real parameter,
$u\in[U-1,U+1]$, and put
$$ E(u):= {2^{\l-2}(\l-1)!\over\rr}\,R(2\pi u^2)\,.
\eqno(5.1) $$ Then according to Zhai\fn{3}{In fact, Zhai in his
notation tacitly assumes that $r_1=\dots=r_\l=1$, which means no
actual loss of generality. We have supplemented the factor $\rr$ in
(5.1).} [21], Lemma 3.1, for $\l$ even,
$$ \eq{E(u) &= E^*(u) + \O{u^{2\l-1}}\,, \cr E^*(u) &:= -
\sum_{1\le m\le u} m(u^2-m^2)^{\l-1}\psi\({u^2\over2m}-{m\over2}\)
\,.\cr} \eqno(5.2)
$$ We apply Lemma 1 in the form $-\psi\ge\Sigma_H - \Sigma_H^*$,
choosing $H=[U]$. Thus we get $$ \eq{E^*(u)\ge &\ - U  + \DS
m(u^2-m^2)^{\l-1}\times\cr &\times\(\alh\sin\(2\pi
h\({u^2\over2m}-{m\over2}\)\)- \beh\cos\(2\pi
h\({u^2\over2m}-{m\over2}\)\)\) \,.\cr} \eqno(5.3)$$ We split up the
range $1\le m\le u$ into dyadic subintervals $\Mj=]M_{j+1},M_j]$,
$M_j=u\,2^{-j}$ for $j=0,\dots,J$, where $J$ is minimal such that
$(U-1)2^{-J-1}<1$. We thus have to deal with exponential sums $$
\E_j(h,u):=\sum_{m\in\Mj}m(u^2-m^2)^{\l-1}\,
e\(-h\({u^2\over2m}-{m\over2}\)\)\,. $$ We transform them by means
of Lemma 2, with $$ G(\xi) = \xi(u^2-\xi^2)^{\l-1}\,,\qquad F(\xi)=
-h\({u^2\over2\xi}-{\xi\over2}\)\,. $$ By straightforward
computations, on each interval $\Mj$ the conditions of Lemma 2 are
fulfilled with the parameters $X= M_j$, $Y={M_j^3\over hu^2}$,
$Z=M_j\,u^{2\ell-2}$. We obtain $$ \eq{\E_j(h,u) = & \
h^{3/4}\,u^{2\l-1/2} \sum_{k\in F'(\Mj)}
{(2k-2h)^{\l-1}\over(2k-h)^{\l+1/4}}\,e\(-u\sqrt{h}\sqrt{2k-h}-{\textstyle{
1\over8}}\)\cr &+ \O{u^{2\l-3}{M_j^{5/2}\over h^{1/2}}+u^{2\l-1}\log
u}\,.\cr}  \eqno (5.4)$$ We first bound the overall contribution of
the error terms, summing over $j$ and $h$. Let $\gah$ denote $\alh$
or $\beh$, thus $\gah\ll h^{-1}$ in any case, then $$ \sum_{j=0}^J
\sum_{h=1}^{[U]} \gah \(u^{2\l-3}{M_j^{5/2}\over
h^{1/2}}+u^{2\l-1}\log u\) \ll u^{2\l-1/2}+u^{2\l-1}(\log u)^3 \ll
u^{2\l-1/2}\,. \eqno (5.5) $$ Summing up the main terms in (5.4), we
notice that the total range of $k$ becomes $h=F'(M_0)< k \le
F'(M_{J+1})=\h h+2^{2J+1} h=:K_{h,U}$, and we obtain $$ \eq{&\DS
m(u^2-m^2)^{\l-1} \gah \,e\(-h\({u^2\over2m}-{m\over2}\)\) = \cr &=
u^{2\l-1/2} \Sh \gah\,h^{3/4}\, \sum_{h<k\le K_{h,U}}
{(2k-2h)^{\l-1}\over(2k-h)^{\l+1/4}}\,e\(-u\sqrt{h}\sqrt{2k-h}-{\textstyle{
1\over8}}\) + \cr & + \O{u^{2\l-1/2}}\,.\cr} \eqno(5.6) $$ Using the
real and imaginary part of this result in (5.3), we arrive at
$$ E^*(u) \ge u^{2\l-1/2}\,S(u,U) - c_1 u^{2\l-1/2}\,, \eqno(5.7) $$
where $$ \eq{S(u,U) :=&\ \sum_{(h,k)\in\D(U)} h^{3/4} \,
{(2k-2h)^{\l-1}\over(2k-h)^{\l+1/4}}\times\cr&\times\(\alh\,
\sin\(2\pi u\sqrt{h}\sqrt{2k-h}+{\textstyle{\pi\over4}}\)-\beh\,
\cos\(2\pi u\sqrt{h}\sqrt{2k-h}+{\textstyle{\pi\over4}}\)\)\,,\cr}
$$ $$ \D(U) := \{(h,k)\in\Z^2:\ 1\le h\le U\,,\ h<k\le
K_{h,U}\ \}\,, $$ and $c_1$ is an appropriate positive constant.
Our next step is to get rid of "most" of the terms of the last
double sum. To this end we use Lemma 3, multiplying $S(u,U)$ by a
Fej\'er kernel $\F_T(u-U)$, where $T$ is a new large parameter,
and integrating over $U-1\le u\le U+1$. We obtain $$  \eq{& I(T,U)
:= \Int_{U-1}^{U+1} S(u,U)\F_T(u-U)\d u = \Int_{-1}^1
S(U+v,U)\F_T(v)\d v = \cr & = \sum_{(h,k)\in\D(U),\ h(2k-h)\le
T^2} h^{3/4} \,
{(2k-2h)^{\l-1}\over(2k-h)^{\l+1/4}}\(1-{\sqrt{h(2k-h)}\over T
}\)\times\cr&\times\(\alh\, \sin\(2\pi
U\sqrt{h(2k-h)}+{\textstyle{\pi\over4}}\)-\beh\, \cos\(2\pi
U\sqrt{h(2k-h)}+{\textstyle{\pi\over4}}\)\)\cr & +
\O{\sum_{(h,k)\in\D(U)} h^{-3/4}
\,{(2k-2h)^{\l-1}\over(2k-h)^{\l+3/4}}}\,.\cr } \eqno(5.8) $$ The
$O$-term here is harmless: In fact, $$ \sum_{h=1}^\infty
h^{-3/4}\sum_{k>h}{(2k-2h)^{\l-1}\over(2k-h)^{\l+3/4}} \ll
\sum_{h=1}^\infty h^{-3/4}\sum_{m>h/2}m^{-7/4} \ll
\sum_{h=1}^\infty h^{-3/2}\ll1\,.\eqno(5.9) $$

We proceed to derive a lower bound for the main term in (5.8). To
this end we relate the two parameters $U$ and $T$ to each other:
For arbitrary $T$ sufficiently large, we choose $U$ according to
the Dirichlet approximation theorem, such that $$ T^2\le U\le
T^2\,16^{T^2} \leqno{\hbox{(i)}} $$ and
$$ \norm{U\sqrt{h(2k-h)}}\le {1\over16} \leqno{\hbox{(ii)}} $$ for
all $(h,k)\in\D(U),\ h(2k-h)\le T^2$, where $\norm{\cdot}$ denotes the
distance from the nearest integer. Now $h(2k-h)\le T^2\,,\ k>h$,
together with (i) implies that
$$ h\le T\le \sqrt{U}\,. \eqno(5.10) $$ By the definitions in Lemma 1,
for all of the $h$ occurring in (5.8), $\alh \asymp h^{-1}$, and
$\beh\asymp U^{-1}$, hence $$ \eq{\alh\,& \sin\(2\pi
U\sqrt{h(2k-h)}+{\textstyle{\pi\over4}}\)-\beh\, \cos\(2\pi
U\sqrt{h(2k-h)}+{\textstyle{\pi\over4}}\)\ge \cr &\ge c_1
\sin\({\pi\over8}\){1\over h}-{c_2\over U}\ge {c_3\over h}\,,\cr}
$$ with suitable positive constants $c_1, c_2, c_3$.
Putting for short $$ \theta_{U,T,\l}(n) := \sum_{(h,k)\in\D(U),\
h(2k-h) = n} {h^{1/2}\over(2k-h)^{1/2}}\,\(1-{h\over
2k-h}\)^{\ell-1}\,,
$$ we thus readily infer from (5.8) and (5.9) that, for some
$c_4, c_5>0$,
$$ \eq{I(T,U) &\ge c_3\sum_{1\le n\le T^2} {\theta_{U,T,\l}(n)\over n^{3/4}}
\(1-{\sqrt{n}\over T}\)\ - c_4\cr &\ge c_5 \sum_{1\le n\le T^2/2}
{\theta_{U,T,\l}(n)\over n^{3/4}} \ - c_4\,.  \cr}\eqno(5.11) $$
We notice further that $(h,k)\in\D(U)$ explicitly means that $$
1\le h\le U\,,\quad h<k\le \h h+2^{2J+1} h \asymp U^2 h\,, $$
while $h(2k-h)\le T^2\,,\ k>h$ implies (5.10) and $$ k<2k-h\le
{T^2\over h} \le T^2\le U\,. $$ Hence, for $1\le n\le T^2$, $$
\eq{&\theta_{U,T,\l}(n) = \theta_\l(n):= \sum_{h(2k-h)=n\atop k>h}
{h^{1/2}\over(2k-h)^{1/2}}\,\(1-{h\over 2k-h}\)^{\ell-1}\cr &\gg
\sum_{h(2k-h)=n\atop k>2h} {h^{1/2}\over(2k-h)^{1/2}} \gg
\sum_{hm=n,\ m>3h\atop h\equiv m \mod2}{\sqrt{h}\over\sqrt{m}}\gg
\sum_{hm=n,\ 3h<m\le4h\atop h\equiv m \mod2}1\,. \cr }   $$
Therefore, by (5.11), $$ \eq{I(T,U) + c_4 &\gg \sum_{T^2/4\le n\le
T^2/2 \atop n\odd} n^{-3/4} \sum_{hm=n\atop 3h<m\le4h } 1 \cr
&\gg\ T^{-3/2} \sum_{T^2/4\le hm\le T^2/2\atop 3h<m\le4h,\ h,m
\odd}1 \ \gg\ T^{1/2}\,.\cr} $$ Thus, $$ I(T,U) \gg T^{1/2}\gg
(\log U)^{1/4}\,,
$$ since (i) readily implies that $T\gg(\log U)^{1/2}$. On the
other hand, it follows from the definition of $I(T,U)$ that $$
I(T,U) \le \(\sup_{U-1\le u\le U+1 }S(u,U)\)\Int_{-1}^1\F_T(v)\d
v\,. $$ Since $$ \Int_{-1}^1\F_T(v)\d v = \Int_{-T}^T\({\sin(\pi
v)\over\pi v}\)^2\d v\le1\,,   $$ this implies that there exists a
value $u^* \in [U-1,U+1]$ for which $$ S(u^*,U)\ge c_6 (\log
u^*)^{1/4}\,, \eqno(5.12) $$ $c_6$ a suitable positive constant.
It remains to recall that if $T$ runs through an unbounded
sequence of positive reals, by construction so do $U$ and $u^*$.
Therefore, (5.1), (5.2), (5.7) and (5.12) together complete the
proof of our theorem. \bsk\msk

\vbox{{\bf 6.~Concluding remarks. } 1. It is appropriate to comment
on the somewhat disturbing restriction that $\l$ must be even. In
general, in (5.2) the argument of the function $\psi$ contains an
additional term $-{\l\over2}$, according to Zhai's [21] Lemma 3.1.
This gives an additional factor $e(\h h\l)$ at the right-hand side
of (5.6), and as a consequence, additional terms $\pi h\l$ in the
arguments of the sine and cosine in (5.8). For $\l$ odd, the
definitions of $\theta_{U,T,\l}(n)$ and $\theta_\l(n)$ ultimately
contain alternating factors $(-1)^h$ which fatally affect our
argument, at least in its present form.}\ssk 2. It is a natural
question whether the sophisticated new methods due to Hafner [4] and
Soundararajan [17] can be applied, in order to improve the result by
a loglog-factor, as it is the case for the divisor and circle
problems. However, both arguments are based on the fact that for
$d(n)$ and $r(n)$, the average order is accomplished by a "thin" set
of integers on which $d(n)$, resp., $r(n)$ attain exceptionally
large values. The improvement is effected by restricting the
application of the Dirichlet approximation theorem to such a thin
set. In our present case, a similar observation concerning the
arithmetic function $\theta_\l(n)$ is at least not at all
straightforward, {\it inter alia } because $\theta_\l(n)$ fails to
be multiplicative.

\bsk  \bsk

{\klein  \parindent=0pt \def\smc{}

\cen{\bf References}  \bsk \hsize=17true cm \vsize=24.5true cm

[1] {\smc F.~Chamizo \and H.~Iwaniec,} On the sphere problem.
Rev.~Mat.~Iberoamericana {\bf 11}, 417-429 (1995). \ssk

[2] C.S.~Gordon and E.N.~Wilson, The spectrum of the Laplacian on
Riemannian Heisenberg manifolds, Michigan Math. J. {\bf33},
253-271 (1986). \ssk

[3] S.W.~Graham and G.~Kolesnik, Van der Corput's method on
exponential sums, Cambridge 1991. \ssk

[4] J.L. Hafner, New omega results for two classical lattice point
problems, Invent. Math. {\bf63}, 181-186 (1981).\ssk

[5] D.R.~Heath-Brown, Lattice points in the sphere, In: Number
theory in progress, Proc. Number Theory Conf. Zakopane 1997, eds.
K.~Gy\"ory et al., vol. 2 (1999), 883-892. \ssk

[6] L.~H\"{o}rmander, The spectral function of an elliptic operator,
Acta Math. {\bf121}, 193-218 (1968).\ssk

[7] {\smc M.N.~Huxley}, {Area, lattice points, and exponential
sums.} LMS Monographs, New Ser. {\bf 13}, Oxford 1996.  \ssk

[8] M.N.~Huxley, Exponential sums and lattice points III.
Proc.~London Math.~Soc. (3) {\bf87}, 591-609 (2003). \ssk

[9] {\smc A.~Ivi\'c, E.~Kr\"{a}tzel, M.~K\"{u}hleitner, \and W.G.~Nowak,}
Lattice points in large regions and related arithmetic functions:
Recent developments in a very classic topic. Proceedings Conf.~on
Elementary and Analytic Number Theory ELAZ'04, held in Mainz, May
24-28, W.~Schwarz and J.~Steuding eds., Franz Steiner Verlag 2006,
pp. 89-128.\ssk

[10] {\smc H.~Iwaniec, E.~Kowalski,} Analytic Number Theory, AMS
Coll.Publ.~53. Providence, R.I., 2004. \ssk

[11] M.~Khosravi, Spectral statistics for Heisenberg manifolds,
Ph.D.~thesis, McGill U. 2005. \ssk

[12] M.~Khosravi and Y.N.~Petridis, The remainder in Weyl's law
for $n$-dimensional Heisenberg manifolds, Proc. AMS {\bf133}/12,
3561-3571 (2005).\ssk

[13] M.~Khosravi and J.A.~Toth, Cramer's formula for Heisenberg
manifolds, Ann.~de l'institut Fourier {\bf55}, 2489-2520
(2005).\ssk

[14] {\smc E.~Kr\"atzel,} Lattice points. Berlin 1988. \ssk

[15] {\smc E.~Kr\"atzel,} Analytische Funktionen in der
Zahlentheorie. Stuttgart-Leipzig-Wiesbaden 2000.  \ssk

[16] Y.N.~Petridis and J.A.~Toth, The remainder in Weyl's law for
Heisenberg manifolds, J.~Diff.~Geom. {\bf60}, 455-483 (2002).\ssk

[17] {\smc K.~Soundararajan,} Omega results for the divisor and
circle problems. Int.~Math.~Res.~Not. {\bf36}, 1987-1998 (2003).
\ssk

[18] {\smc K.-M.~Tsang,} Counting lattice points in the sphere.
Bull.~London Math.~Soc. {\bf32}, 679-688 (2000). \ssk

[19] J.D.~Vaaler, Some extremal problems in Fourier analysis,
Bull.~Amer.~Math.~Soc. {\bf12}, 183-216 (1985).\ssk

[20] A.~Walfisz, Gitterpunkte in mehrdimensionalen Kugeln, Warszaw
1957.\ssk

[21] W.~Zhai, On the error term in Weyl's law for the Heisenberg
manifolds, Acta Arithm.~{\bf134}, 219-257 (2008).

\vbox{\vskip 0.7true cm}

\parindent=1.5true cm

\vbox{Institute of Mathematics

Department of Integrative Biology

Universit\"at f\"ur Bodenkultur Wien

Gregor Mendel-Stra\ss e 33

1180 Wien, \"Osterreich \ssk

E-mail: {\tt nowak@boku.ac.at} \ssk

}}

\bye